\documentclass{amsart}
\usepackage{amssymb,latexsym}
\usepackage{newlattice}
\usepackage{graphicx}

\renewcommand{\Ji}{J}

\newtheorem{theorem}{Theorem}

\theoremstyle{definition} 

\theoremstyle{plain}

\begin{document}
\title[Principal congruences of a finite lattice]
{Revisiting the representation theorem of\\
finite distributive lattices\\
with principal congruences.\\
A \emph{Proof-by-Picture} approach}  
\author{G. Gr\"{a}tzer} 
\email[G. Gr\"atzer]{gratzer@me.com}
\urladdr[G. Gr\"atzer]{http://server.maths.umanitoba.ca/homepages/gratzer/}

\author[H. Lakser]{H. Lakser}
\email[H. Lakser]{hlakser@gmail.com} 
\address{Department of Mathematics\\University of Manitoba\\Winnipeg, MB R3T 2N2\\Canada}

\date{May 9, 2020}

\begin{abstract}
A classical result of R.\,P. Dilworth states that every finite distributive lattice $D$  
can be represented as the congruence lattice of a finite lattice~$L$.
A~sharper form was published in G.~Gr\"atzer and E.\,T. Schmidt in 1962,
adding the requirement that all congruences in $L$ be principal.
Another variant, published in 1998 by the authors and E.\,T. Schmidt, 
constructs a planar semimodular lattice $L$.
In this paper, we merge these two results:
we construct $L$ as a planar semimodular lattice in which all congruences are principal.
This paper relies on the techniques developed by the authors and E.\,T. Schmidt in the 1998 paper.
\end{abstract}

\maketitle

\section{Introduction}\label{S:Introduction}

Let us start with the classical result 
of R.\,P. Dilworth from 1942 (see the book \cite{BFK90}
for background information):

\begin{theorem}\label{T:Dilworth}
Every finite distributive lattice $D$ 
can be represented as the congruence lattice 
of a finite lattice $L$.
\end{theorem}

A sharper form was published in 
G.~Gr\"atzer and E.\,T. Schmidt \cite{GS62} 
(see also Theorem 8.5 in \cite{CFL2}).
The new idea was the use of standard ideals, 
see G.~Gr\"atzer~\cite{GS59a}
and G. Gr\"atzer and E.\,T. Schmidt~\cite{GS61}.

\begin{theorem}\label{T:full}
Every finite distributive lattice $D$ 
can be represented as the congruence lattice 
of a finite relatively complemented lattice $L$.
\end{theorem}

All congruences are principal 
in a finite relatively complemented lattice $L$. 
So~we obtain the following variant of Theorem~\ref{T:full}.

\begin{theorem}\label{T:principal}
Every finite distributive lattice $D$ 
can be represented as the congruence lattice 
of a finite lattice $L$ in which all congruences are principal.
\end{theorem}

G. Gr\"atzer, H.~Lakser, 
and E.\,T. Schmidt \cite{GLS98a}
proved another variant of Theorem~\ref{T:full}.

\begin{theorem}\label{T:planar}
Let $D$ be a finite distributive lattice.  
Then there exists a planar semimodular lattice $L$  
with $\Con L$ isomorphic to $D$. 
\end{theorem}

In this note, we combine Theorem~\ref{T:principal}
and~\ref{T:planar},
using the techniques developed for Theorem~\ref{T:planar}.

\begin{theorem}\label{T:main}
Every finite distributive lattice $D$ 
can be represented as the congruence lattice 
of a planar semimodular lattice $L$
in which all congruences are principal.  
\end{theorem}

There are other aspects of these constructions
discussed in the book \cite{CFL2}, 
for instance, the size of $L$. 
The constructions in Theorems~\ref{T:Dilworth},
\ref{T:full}, and~\ref{T:main} are ``large'' (exponential),
in Theorem~\ref{T:planar} they are small (cubic polynomial).

There are related results in 
G. Gr\"atzer and H. Lakser~\cite{GL18} and \cite{GL19}.

\subsection*{Outline}
For a formal proof of Theorem~\ref{T:main}, 
we need the formal proof of Theorem~\ref{T:planar},
as presented in G. Gr\"atzer, H.~Lakser, 
and E.\,T. Schmidt \cite{GLS98a}.
There are two obvious solutions: 
copy the formal proof from \cite{GLS98a}
(making the editor unhappy)
or require that the reader be familiar
with the paper \cite{GLS98a}
(making the reader unhappy).
So we choose the middle ground, 
we present a \emph{Proof-by-Picture} 
(as defined in \cite{CFL2}) of Theorem~\ref{T:planar}.
We do this in Section~\ref{S:pbp} and
complete the proof of Theorem~\ref{T:main} 
in~Section~\ref{S:Th5}.

\subsection*{Notation}
We use the notation as in \cite{CFL2}.

In particular, for the ordered sets $P$ and $Q$,
we can form the (ordinal) sum, $P + Q$ 
and the glued sum $P \dotplus Q$,
as illustrated in Figure~\ref{F:PandQ}.
Observe that the glued sum $P \dotplus Q$ requires 
that $P$ has a unit and $Q$ has a zero (which are identified). 

\emph{Coloring} of a finite lattice $L$ 
attaches a join-irreducible congruence to an edge 
(covering interval) of $L$ generating it,
see Figures~\ref{F:D}--\ref{F:NandS}
for examples.

\section{\emph{Proof-by-Picture} of Theorem~\ref{T:planar}}
\label{S:pbp}

\begin{figure}[p]
\centerline{\includegraphics[scale=1]{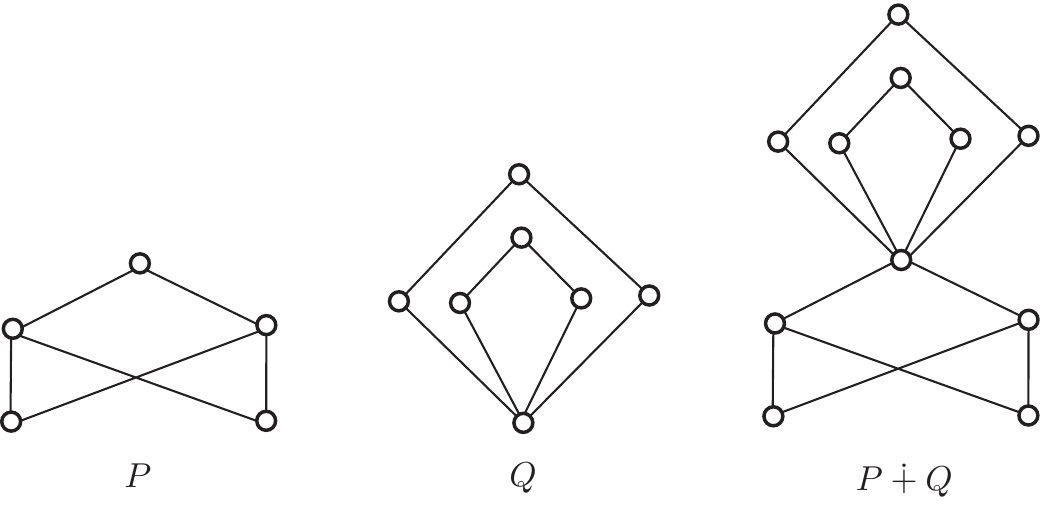}}
\caption{Glued sum of two ordered sets, $P$ and $Q$}
\label{F:PandQ}

\medskip

\medskip

\medskip

\medskip

\includegraphics[scale=1]{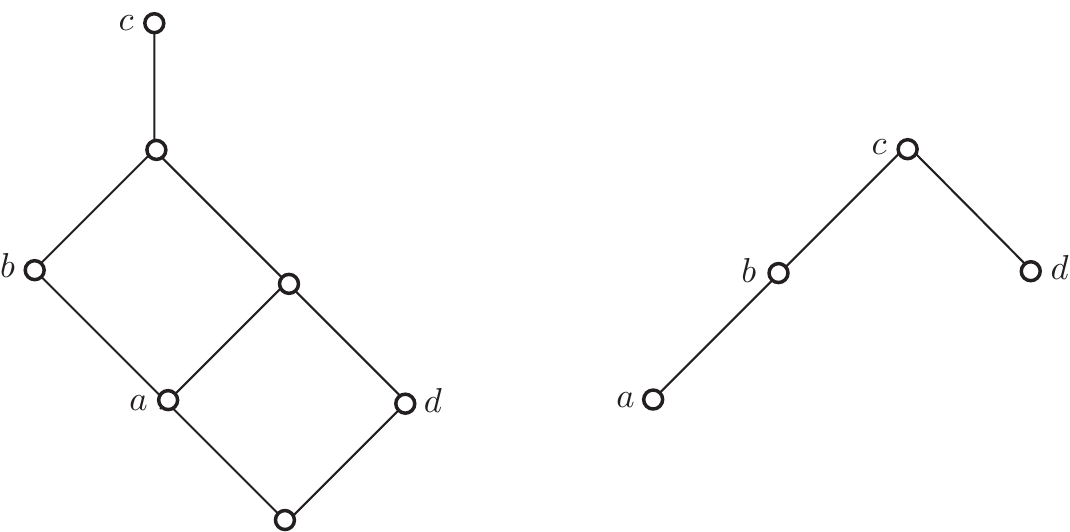}
\caption{The lattice $D$ to represent and 
the ordered set $P = \Ji(D)$}
\label{F:D}

\medskip

\medskip

\medskip

\medskip

\centerline{\includegraphics[scale=1]{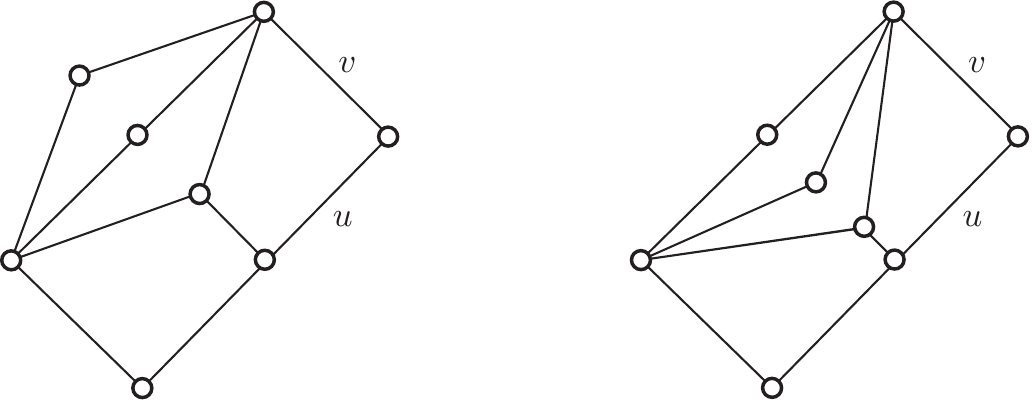}}
\caption{Two diagrams of the building block $S_8(u,v)$, $u \prec v$}
\label{F:S8}
\end{figure}

\begin{figure}[p]
\centerline{\includegraphics[scale = 1.0]{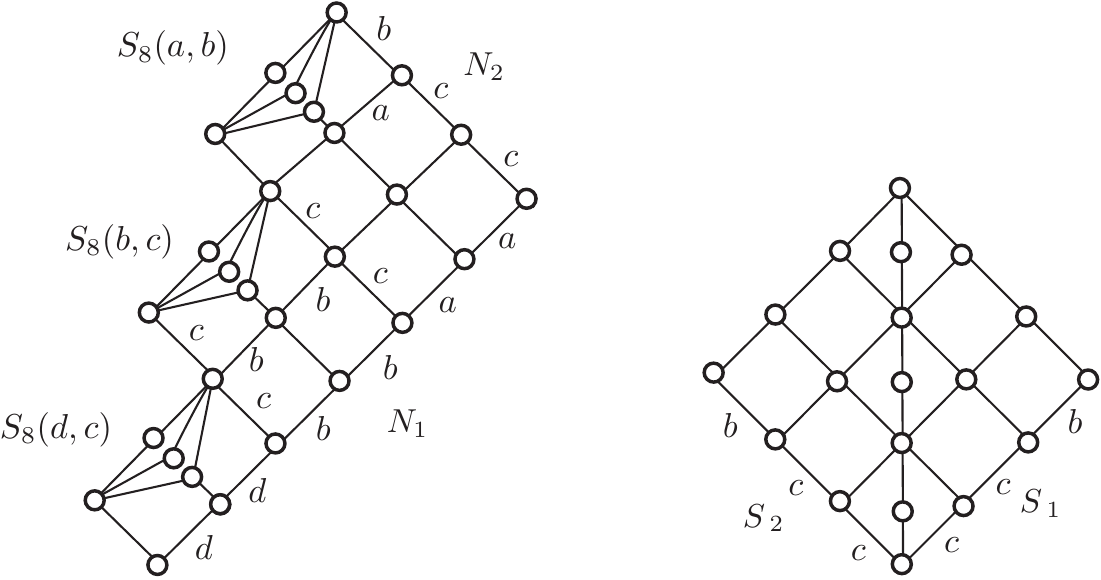}}
\caption{The lattices $N$ and $S$}
\label{F:NandS}

\vspace{30pt}

\centerline{\includegraphics[scale = 0.8]{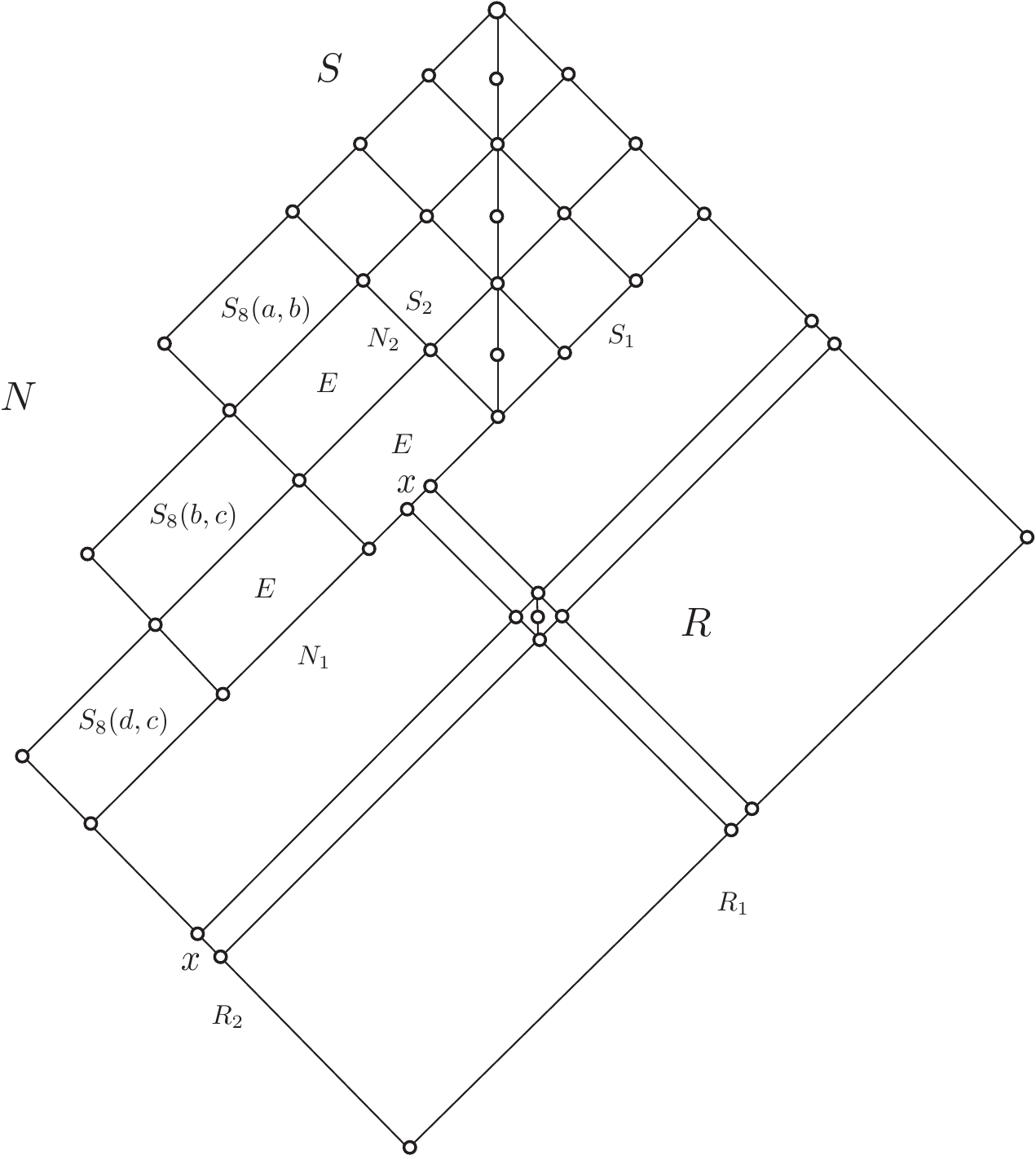}}
\caption{A sketch of $L$ without the ``tail''}
\label{F:Sketch}
\end{figure}

We start constructing the planar semimodular lattice $L$
of Theorem~\ref{T:main} 
for the distributive lattice $D$ and the ordered set $P = \Ji(D)$
of Figure~\ref{F:D}, with the three lattices, 
the planar semimodular lattices $N$ (for Nondistributive), 
$S$ (for Square), and~$R$ (for Rectangle).
We glue them together and add some covering $\SM 3$-s,
to~obtain~$L$, as sketched in Figure~\ref{F:Sketch}.

In Steps 1--4, we assume that $P$ has no \emph{isolated elements},
that is, for every $x \in P$, there is a $y \in P$ 
with $x < y$ or $y < x$.
 
\emph{Step 1: Constructing $N$.}
Take the eight-element, planar, 
semimodular lattice $S_8$ of Figure~\ref{F:S8}. 
We take three copies, $S_8(a,b)$, $S_8(b,c)$, $S_8(d,c)$,
one for every covering pair in $P = \Ji(D)$.
Let $E = \SC 2 \times \SC 3$.
We glue these together (preserving the colors!)
as in Figure~\ref{F:NandS}.
More precisely, we glue $S_8(b,c)$ to $E$, 
and glues $S_8(d,c)$ to the top left boundary of $E$.
Then we glue $D$ to this lattice twice
and glue $S_8(a,b)$ to the top.
We denote by $N_1$ and $N_2$ the lower right and the upper right
boundaries of $N$, respectively.

\emph{Step 2: Constructing $S$.}
We~form $N_2^2$. 
In every covering square of the main vertical diagonal, 
we add an element to make it an~$\SM 3$,
forming the lattice $S$, see Figure~\ref{F:NandS}.  
We denote by $S_1$ and $S_2$ the lower left and lower right
boundaries of $S$, respectively.
This will make a copy of the colors $b$ and $c$ in $S_2$,
making them available for the $\SM 3$ insertions in Step 4b.

\emph{Step 3: Constructing $R$.}
Let the chain $C_1$ be isomorphic to $N_1 \dotplus S_1$.
We choose a chain $C$ of length four 
and color the edges with $\set{a,b,c,d}$ 
(in any order). Define $R = C \times C_1$. 
We denote by $R_1$, $R_2$, and $R_1'$ 
the lower right, lower left, and upper left boundaries of $R$, respectively.

\emph{Step 4: Constructing $L$.}

\emph{Step 4a: Gluing $N$, $S$, and $R$.}
We glue $N$ and $S$ by identifying $N_2$
with $S_2$ (preserving colors!); we call this lattice $L_1$.
Then we glue $L_1$ and $R$ by identifying $R_1'$
with the lower right boundary of $L_1$ (preserving colors!);
let $L_2$ be the lattice we obtain.

\emph{Step 4b: Adding $\SM 3$-s to $L_2$.} 
Every color $x$ occurs in $N_1\dotplus S_1 = R_1'$
as the color of an edge. 
If $x$ is not a maximal element in $P$, 
Then $x$ occurs in $N_1$ as the color of~an edge (maybe many times).
If $x$ is a maximal element in $P$, 
then $x$ occurs in~$S_1$ as the color of an edge (maybe many times),
so $x$ occurs in $S_2$ as the color of an edge,
and therefore also in $R_1'$.

So in the grid $R$, we take a ``covering row'' 
and a ``covering column'' 
hitting $R_1'$ and $R_2$ in edges of color $x$,
see Figure~\ref{F:Sketch}.
They determine a covering square to which we add an element
to obtain an $\SM 3$.
We do this for all covering squares given by 
a covering row and a covering column both colored by $x$, 
thereby identifying all the principal congruences 
determined by a prime interval colored by $x$.

We repeat this for every color $x$.

The $S_8(u,v)$ sublattices then determine 
the desired order on the join-irreducible
congruences---see Figure~\ref{F:S8}.

\emph{Step 5: Adding the tail.} 
If there are $k > 0$ isolated elements,  
we form $\SC {k-1} \dotplus L$; 
the tail is $\SC {k-1}$.

This completes the \emph{Proof-by-Picture} 
of Theorem~\ref{T:planar}.

\section{Proving Theorem~\ref{T:main}}\label{S:Th5}

We have to modify the construction of the
planar semimodular lattice $L$ of Section~\ref{S:pbp}
to make all congruences principal. 
In Step 3, we choose a chain $C$ of length four. 
Observe that the proof of Theorem~\ref{T:planar}
remains valid as long as every color 
is represented as the coloring of $C$.

Now we change the definition of $C$. 
For every $x \in D$, define 
\[
   r(x) = \setm{a \in \Ji(D)}{x \leq a},
\]
and let $C_x$ be a chain of $|r(x)|+1$ elements,
colored by the elements of $r(x)$ (in any order).
Let $0_x, 1_x$ denote the bounds of $C_x$.
Let $C$ be the glued sum of the chains 
$C_x$ for $x \in D$ (in any order). 
This chain $C$ obviously satisfies 
the condition that every color 
is represented as the color of an edge in $C$.

Therefore, the lattice $L$ constructed in Section~\ref{S:pbp}
satisfies the requirements of Theorem~\ref{T:planar}.
We only have to observe that all congruences are principal.  

Let $\bga$ be a congruence of $L$. 
Let $x$ be an element of $D$ that corresponds to $\bga$
under an isomorphism between $\Con L$ and $D$.
Since $C_x$ is colored by the set $r(x)$,
we conclude that in $L$, we have 
\[
   \con{0_x, 1_x} = \bga,
\]
completing the proof.

\end{document}